\documentclass{article}

\usepackage{amsfonts,amsmath,amssymb}
\usepackage{graphicx}
\usepackage{natbib}
\usepackage{bm}

\title{An eddifying Stommel model: Fast eddy effects in a two-box ocean}

\author{
  William Barham
  and
  Ian Grooms\thanks{Corresponding author. Email: ian.grooms@colorado.edu}\\\vspace{6pt}
Department of Applied Mathematics, \\
University of Colorado,\\
Boulder, CO 80309}

\begin{document}

\maketitle

\begin{abstract}
A system of stochastic differential equations is formulated describing the heat and salt content of a two-box ocean. 
Variability in the heat and salt content and in the thermohaline circulation between the boxes is driven by fast Gaussian atmospheric forcing and by ocean-intrinsic, eddy-driven variability. 
The eddy forcing of the slow dynamics takes the form of a colored, non-Gaussian noise.
The qualitative effects of this non-Gaussianity are investigated by comparing to two approximate models: one that includes only the mean eddy effects (the `averaged model'), and one that includes an additional Gaussian white-noise approximation of the eddy effects (the `Gaussian model').
Both of these approximate models are derived using the methods of fast averaging and homogenization. 

In the parameter regime where the dynamics has a single stable equilibrium the averaged model has too little variability.
The Gaussian model has accurate second-order statistics, but incorrect skew and rare-event probabilities.
In the parameter regime where the dynamics has two stable equilibria the eddy noise is much smaller than the atmospheric noise.
The averaged, Gaussian, and non-Gaussian models all have similar stationary distributions, but the jump rates between equilibria are too small for the averaged and Gaussian models.
\end{abstract}

\section{Introduction}
H.~\citet{Stommel61} developed a conceptual model of the global ocean thermohaline circulation that consists of a system of ordinary differential equations modeling the heat and salt content of two containers (`boxes'). 
One box models the equatorial ocean, and the other models the extra-tropical ocean.
The boxes exchange heat and salt with each other and with the atmosphere.
The rate of flow between the boxes is proportional to the density difference between the boxes, and a major result of Stommel's investigation was that in some parameter regimes the system exhibits two equilibria: one analogous to the current climate, with dense cold water sinking at high latitudes, and one corresponding to a very different regime with dense salty water sinking in the equatorial ocean.
In general, the goal of studies using extremely simplified models like Stommel's is to observe and understand qualitative features that might inform and guide subsequent studies using more complete and more complex models.
The qualitative predictions of Stommel's model have since been verified using more complete ocean models, e.g.~\citet{Rahmstorf95} and \citet{DeshayesEtAl13}.

The present investigation develops a model closely related to Stommel's where the slow, density-driven exchange of heat and salt between the boxes is augmented by fast, non-Gaussian stochastic processes representing eddy-driven heat and salt transport.
Eddies smaller than the grid scale of comprehensive numerical ocean (and atmosphere) models can have significant impacts on the global circulation, and modeling the impacts of these unresolved eddies is a topic of continuing research; \citet{BernerEtAl16} contains a review of stochastic models of eddy effects from an operational modeling perspective.
The second author recently proposed a non-Gaussian model of the heat and salt transport associated with unresolved ocean eddies \citep{Grooms16}.
In this model, the eddy velocity and density fields (the latter linearly related to temperature and salinity) are modeled as centered Gaussian random fields, and the transports are modeled as the product of eddy velocity and density.
The product of centered, jointly Gaussian random variables has a distinctive, non-Gaussian probability density, with a logarithmic singularity at the origin and skewed, algebraically-modulated exponential decay in the tails.
This non-Gaussian model is significantly different from recent Gaussian stochastic models of eddy transport, e.g.~\citet{ACJPWZ16,WHGSJ16} and \citet{JPZ17}.
The present investigation is motivated by the desire to observe the qualitative effects of this kind of non-Gaussian transport in an extremely simple model, in particular by comparison to Gaussian stochastic models, with the expectation of informing future investigations using more complex models.

Several authors have developed stochastic versions of Stommel's model to investigate the slow response of the ocean thermohaline circulation to fast atmospheric forcing, e.g.~\cite{Cessi94,VACTH01,Monahan02,MTL02} and \citet{MC11}.
In these stochastic Stommel models the atmospheric heat and salt fluxes in Stommel's model are replaced by Gaussian stochastic noise terms, resulting in a system of stochastic differential equations (SDEs).
The model developed here attempts to understand a qualitatively different physical process: fast eddy transport.
Since the eddies are typically faster than the global thermohaline circulation, the new model has the form of a slow-fast system, where eddy variables evolve on a fast time scale and converge towards a jointly Gaussian distribution conditioned on the slow variables.
The slow variables (the heat and salt difference between the boxes) are impacted by quadratic products of fast variables modeling the fast eddy transport.
The formal theory of fast averaging \citep{PK74,PS08,FW12}, is used to generate approximate slow systems for comparison: one with a drift correction and one with both drift and diffusion corrections derived from the eddy dynamics.
These approximate systems qualitatively represent more complete ocean models with, respectively, deterministic and Gaussian stochastic models of the eddy transport.

A new stochastic Stommel model is developed in \S \ref{sec:Model}.
The two approximate models of the slow system are derived in \S \ref{sec:Approx}.
The numerical methods and experimental configuration are described in \S \ref{sec:methods} and the results of these simulations are described in \S \ref{sec:results}.
A slightly different model with two stable equilibria is formulated and simulated in \S\ref{sec:NoMeanDiffusion}.
The results and their implications are discussed in \S \ref{sec:conclusions}.

\section{Formulating a Slow-Fast Two-Box Stochastic Ocean Model}
\label{sec:Model}
Consider a domain $[0,L_x]\times[0,L_y]\times[0,H]$ representing an ocean basin, and let this domain be partitioned into two subdomains $[0,L_x]\times[0,\ell]\times[0,H]$ and $[0,L_x]\times[(\ell,L_y]\times[0,H]$ with volumes $V_1=L_x\ell H$ and $V_2=L_x(L_y-\ell)H$.
The first box (index 1) will represent the equatorial side of the ocean basin, and the second will represent the poleward side.
The domain is filled with a fluid whose density is related to its temperature and salinity via
\[\rho = \rho_0[1 + \alpha_S(S-S_0) -\alpha_T(T-T_0)]\]
where $\rho_0=1029$ kg/m$^3$ is a constant reference density, $T_0=5$ C and $S_0=35$ psu are a constant reference temperature and salinity (psu are practical salinity units; for the present purposes it is reasonable to use the simplification 1 psu = 1 g/kg), and $\alpha_S=7.5\times10^{-4}$ psu$^{-1}$ and $\alpha_T=1.7\times10^{-4}$ C$^{-1}$ are coefficients of haline and thermal expansion.
The conservation equations for heat are in the form of a system of two differential equations
\[\frac{\text{d}T_1}{\text{d}t} = -\frac{1}{\tau_T}(T_1-T_1^*) - \frac{F_{q}}{\rho_0c_pV_1},\;\;\frac{\text{d}T_2}{\text{d}t} = -\frac{1}{\tau_T}(T_2-T_2^*) + \frac{F_{q}}{\rho_0c_pV_2}\]
where $T_1$ and $T_2$ are the mean temperature in each box, $\tau_T$ is the timescale of relaxation towards an externally-specified atmospheric temperature $T_i^*$, $c_p=4000$ J/kg is the heat capacity of seawater (e.g.~$\rho_0 c_p V_1 T_1$ is the heat content of the equatorial box), and $F_q$ is the heat flux from the equatorial box to the poleward box.
The total heat content $\rho_0c_p(V_1T_1+V_2T_2)$ thus depends only on the external forcing.

Similarly, the conservation equations for salt are
\[\frac{\text{d}S_1}{\text{d}t} = \frac{1}{2}F(t) - F_{S},\;\;\frac{\text{d}S_2}{\text{d}t} = -\frac{1}{2}F(t) + F_{S}\]
where $2F(t)$ is the external freshwater forcing in the equatorial box (e.g.~rain, runoff, ice melt) and $F_S$ is the salt flux from the equatorial box to the poleward box.
The external freshwater forcing is assumed not to change the net salt content, so that $S_1+S_2$ remains constant in time.

Following \citet{Stommel61}, the heat and salt fluxes between the boxes are assumed to depend only on the temperature and salinity differences between the boxes.
As a result, the temperature and salinity differences between the boxes decouple from the net heat and salt content.
Defining $\Delta T = T_1-T_2$ and $\Delta S = S_1-S_2$, 
\[\frac{\text{d}\Delta T}{\text{d}t} = -\frac{1}{\tau_T}(\Delta T - \Delta T^*) - 2\left[\frac{1}{\rho_0c_pV_1}+\frac{1}{\rho_0c_pV_2}\right]F_T\]
\[\frac{\text{d}\Delta S}{\text{d}t} = F(t)-2F_S.\]
Similar to \citet{Cessi94} and \citet{VACTH01}, the atmospheric temperature difference $\Delta T^*$ and external freshwater forcing are here modeled as constant mean terms plus Gaussian white noise, leading to
\[\text{d}\Delta T= \left[-\frac{1}{\tau_T}(\Delta T - \overline{\Delta T^*}) - 2\left[\frac{1}{\rho_0c_pV_1}+\frac{1}{\rho_0c_pV_2}\right]F_T\right]\text{d}t + \frac{\sigma_{\Delta T}}{\sqrt{\tau_T}}\text{d}W_{\Delta T}\]
\[\text{d}\Delta S = \left[\overline{F}-2F_S\right]\text{d}t+\frac{\sigma_S}{\sqrt{\tau_d}}\text{d}W_S.\]
The amplitude of the noise forcing is here related to $\tau_d$, a diffusive time scale, defined below.

In Stommel's original model the fluxes between the boxes consists of diffusive fluxes proportional to the temperature and salinity differences, and advective fluxes associated with the large-scale ocean circulation whose rate is proportional to the magnitude of the density difference between the boxes
\[2\left[\frac{1}{\rho_0c_pV_1}+\frac{1}{\rho_0c_pV_2}\right]F_T =\left(\frac{1}{\tau_d}+\frac{1}{\tau_a\rho_0\alpha_T\overline{\Delta T^*}}|\Delta\rho|\right)\Delta T\]
\[2F_S =\left(\frac{1}{\tau_d}+\frac{1}{\tau_a\rho_0\alpha_T\overline{\Delta T^*}}|\Delta\rho|\right)\Delta S\]
where $\tau_d$ is the time scale of diffusive transport, $\tau_a$ is the time scale of advective transport, and
\[\Delta\rho = \rho_0[\alpha_S\Delta S-\alpha_T\Delta T]\]
is the density difference between the boxes.
\citet{Cessi94} used a smoother formulation, which does not qualitatively change the results
\[2\left[\frac{1}{\rho_0c_pV_1}+\frac{1}{\rho_0c_pV_2}\right]F_T =\left(\frac{1}{\tau_d}+\frac{1}{\tau_a(\rho_0\alpha_T\overline{\Delta T^*})^2}\Delta\rho^2\right)\Delta T\]
\[2F_S =\left(\frac{1}{\tau_d}+\frac{1}{\tau_a(\rho_0\alpha_T\overline{\Delta T^*})^2}\Delta\rho^2\right)\Delta S.\]

The novel contribution to the model made here consists of the addition of fast variables crudely representing eddy velocity $v$, temperature $T$, and salinity $S$ anomalies at the interface between the boxes.
The eddy-induced fluxes between the boxes will be modeled as an addition to the slow diffusive and advective fluxes
\[2\left[\frac{1}{\rho_0c_pV_1}+\frac{1}{\rho_0c_pV_2}\right]F_T =\left(\frac{1}{\tau_d}+\frac{1}{\tau_a(\rho_0\alpha_T\Delta T^*)^2}\Delta\rho^2\right)\Delta T + \left[\frac{1}{\ell}+\frac{1}{L_y-\ell}\right]vT\]
\[2F_S =\left(\frac{1}{\tau_d}+\frac{1}{\tau_a(\rho_0\alpha_T\Delta T^*)^2}\Delta\rho^2\right)\Delta S+\left[\frac{1}{\ell}+\frac{1}{L_y-\ell}\right]vS.\]
The prefactors of $\ell^{-1}+(L_y-\ell)^{-1}$ account for the fact that the boxes need not have equal volume, and that total heat and salt need to be conserved.
For simplicity, only $\ell=L_y/2$ is considered from here on.

In general the flux between the boxes should be described by $\int_0^{L_x}\int_0^HvT$d$z$d$x$ where $v$ and $T$ are evaluated at $y=\ell$.
Our formulation amounts to a severe simplification that ignores the spatial structure of the eddy velocity and temperature perturbations between the boxes, and considers them only as zero-mean jointly-Gaussian variables.
This level of simplification is consistent with the simplification of the ocean to two well-mixed boxes in the original Stommel model, and is guided by the desire to investigate the qualitative effects of Gaussian-product noise, since eddy noise with this structure was recently observed by \citet{Grooms16}.

The fast eddy velocity will be modeled as an Ornstein-Uhlenbeck process
\[\text{d}v = -\frac{1}{\tau_e}v\text{d}t + \sqrt{\frac{2}{\tau_e}}\sigma_v\text{d}W_v\]
where $\tau_e$ is the eddy time scale and $\sigma_v$ is the eddy velocity scale, chosen to be $15$ days and $10$ cm/s, respectively \cite{Stammer97}.
The eddy velocity can be thought of as being set by wind-driven processes independent of the density difference between the boxes.
This is a simplification of the more complex reality where eddy kinetic energy depends also on the large-scale density gradient.
The following model of the eddy dynamics is perhaps more qualitatively appropriate
\begin{align*}
\text{d}v &=-\frac{v}{\tau_e}\text{d}v+\sqrt{\frac{2(1+\mu\Delta\rho^2)}{\tau_e}}\text{d}W_v
\end{align*}
where $\mu>0$ is a parameter representing the sensitivity of the eddy variance to the large-scale density gradient.
This model is not pursued further here, in part because of the difficulties in guaranteeing its ergodicity and in finding a robust numerical method for its solution.

The eddy temperature and salinity anomalies will be modeled as resulting from eddy transport across the large-scale gradients
\[\frac{\text{d}T}{\text{d}t} = -\frac{T}{\tau_e}-v\frac{2\Delta T}{L_y},\]
\[\frac{\text{d}S}{\text{d}t} = -\frac{S}{\tau_e}-v\frac{2\Delta S}{L_y}.\]
The relaxation towards zero on a time scale of $\tau_e$ qualitatively represents the full range of dissipative processes acting on temperature and salinity anomalies: cascade towards small scales and eventual diffusion, and atmospheric damping of thermal anomalies, etc.
The time scale $\tau_e$ should not be associated with any particular physical process, but instead guarantees decorrelation of eddy anomalies on the time scale $\tau_e$.

The governing equations are nondimensionalized using the diffusive time scale, the external constant atmospheric temperature difference $\overline{\Delta T^*}\approx20$ C for large-scale temperature, and the convenient salinity scale $\alpha_T\overline{\Delta T^*}/\alpha_S\approx$ 4.5 psu for large-scale salinity.
It will be convenient to scale the eddy temperature and salinity variables differently; specifically, $T$ will have dimensions $\overline{\Delta T^*}L_y/(\sigma_v\tau_d)$ and $S$ will have dimensions $\alpha_T\overline{\Delta T^*}L_y/(\alpha_S\sigma_v\tau_d)$.
The reason for this unexpected scaling will be commented on shortly.

Following traditional notation, the nondimensional temperature difference will be denoted $x$ and the nondimensional salt difference will be denoted $y$; risking confusion, the nondimensional eddy variables will have the same notation as their dimensional counterparts.
The complete nondimensional system is therefore
\begin{subequations}\label{eqn:True}
\begin{align}
\text{d}x &= \left[-\frac{1}{\epsilon_T}(x-1)-[1+P_a(x-y)^2]x + 4vT\right]\text{d}t + \sqrt{\frac{1}{\epsilon_T}}\sigma_x\text{d}W_x\label{eqn:x}\\
\text{d}y &= \left[1-[1+P_a(x-y)^2]y + 4vS\right]\text{d}t +\sigma_y\text{d}W_y\label{eqn:y}\\
\text{d}v &=-\frac{v}{\epsilon}\text{d}v+\sqrt{\frac{2}{\epsilon}}\text{d}W_v\label{eqn:v}\\
\text{d}T&=-\frac{1}{\epsilon}\left[T+2P^2vx\right]\text{d}t\label{eqn:T}\\
\text{d}S&=-\frac{1}{\epsilon}\left[S+2P^2vy\right]\text{d}t\label{eqn:S}
\end{align}
\end{subequations}
where
\[\epsilon_T = \frac{\tau_T}{\tau_d},\;\;\epsilon = \frac{\tau_e}{\tau_d},\;\;P_a=\frac{\tau_d}{\tau_a},\;\;P_e = \frac{\sigma_v\tau_d}{L_y},\;\;P = \sqrt{\epsilon}P_e.\]
$P_a$ and $P_e$ are P\'eclet numbers comparing the time scales of large-scale advective transport and fast eddy transport to the time scale of diffusion, respectively.

The following parameter estimates are drawn from \cite{Cessi94} and \cite{VACTH01}, and are consistent with the more recent observational analysis of \cite{Schmitt08}.
The diffusive time scale $\tau_d$ is approximately 219 years, and the time scale of large scale advection $\tau_a$ is approximately 35 years.
Cessi estimates $\tau_T$ to be 25 days \cite{Cessi94}, but V\'elez-Belchi et al.~argue convincingly that large-scale temperature anomalies are damped on a slower time scale of approximately $220$ days.
V\'elez-Belchi et al.~used salinity noise whose nondimensional amplitude is here $\sigma_y=0.15$, and assuming that fast atmospheric temperature fluctuations lead to perturbations on the order of 0.07 C implies nondimensional thermal noise has amplitude $\sigma_x=0.005$.
Finally, using a length scale appropriate to the global oceans $L_y\approx 8,250$ km leads to the following set of parameters which are adopted for the remainder of the investigation
\begin{equation}\label{eqn:parameters}
\epsilon_T = \frac{1}{400},\;\;\epsilon = \frac{1}{5000},\;\;P_a = 6,\;\;P_e = 80,\;\;\sigma_x=0.005,\;\;\sigma_y=0.15.
\end{equation}
The reason for scaling $S$ and $T$ differently from $\Delta S$ and $\Delta T$ should now be clear: $2P_e^2$ is the same order of magnitude as $\epsilon^{-1}$, implying that both terms in the evolution equations for $S$ and $T$ are of comparable magnitude.

For the parameters (\ref{eqn:parameters}) the system (\ref{eqn:True}) has three equilibria, two of which are stable.
The equilibria all have $v,T,S=0$, and the stable equilibria occur at $(x,y)\approx(0.989,0.22)$ and $(x,y)\approx(0.998,1.00)$.
In the absence of eddy dynamics, one would expect small atmospheric noise to lead to jumping between the two stable equilibria of the system; this was the focus of \cite{Cessi94,Monahan02,MTL02} and \cite{MC11}.
The existence of stable equilibria is intrinsically tied to the nonlinear terms that model slow advective exchange between the boxes.
As the exchange between the boxes becomes dominated by diffusion instead of advection ($P_a\to0$) one of the stable equilibria disappears in a reverse saddle-node bifurcation leaving a single stable equilibrium.\\

Equations (\ref{eqn:T}) and (\ref{eqn:S}) lack noise terms, implying that the classical conditions for ergodicity \citep{Khasminskii12} do not apply.
Conditions for ergodicity of this type of system of SDEs can be found in \cite{MSH02}.
The first condition is that there is an inner-product norm $\|\bm{\cdot}\|$ such that $\langle\bm{u},\bm{F}(\bm{u})\rangle\le\alpha-\beta\|\bm{u}\|^2$ for some $\alpha,\beta>0$ where $\bm{u}$ is a vector containing the dependent variables and $\bm{F}(\bm{u})$ is the drift.
It is straightforward to verify that $\|\bm{u}\|^2 = x^2+y^2+\epsilon v^2+(2\epsilon/P^2)(T^2+S^2)$ satisfies this condition.
The second condition is that the vectors $\left\{\bm{\rho}_i,[[\bm{F},\bm{\rho}_j],\bm{\rho}_k]\right\}$ span $\mathbb{R}^5$ where $\bm{\rho}_i$, $i=1,\,2,\,3$ are the columns of the diffusion matrix, which are here proportional to the first three standard basis vectors, and $[\mathbf{\cdot},\mathbf{\cdot}]$ is a Lie bracket.
Since $[[\bm{F},\bm{\rho}_1],\bm{\rho}_3]$ and $[[\bm{F},\bm{\rho}_2],\bm{\rho}_3]$ are proportional to the fourth and fifth standard basis vectors, respectively, the system satisfies the conditions of \cite{MSH02} for ergodicity.

\section{Two Approximate Slow Models}
\label{sec:Approx}
In this section two systems of SDEs are derived approximating the evolution of the slow variables $x$ and $y$ in (\ref{eqn:True}).
The system of SDEs (\ref{eqn:True}) with parameters (\ref{eqn:parameters}) has three time scales since $\epsilon<\epsilon_T\ll1$: $x$ evolves significantly more quickly than $y$, yet slower than the eddy variables $v$, $T$, and $S$.
Many previous investigations (which lacked the eddy variables) accounted for the scale separation somewhat crudely by setting $x=1$, and focused on the dynamics of the slowest variable $y$, e.g.~\cite{Cessi94,Monahan02,MTL02}, and \citet{MAW08}.
The analysis of \cite{MC11} is more careful, employing the same methods used here but for the system without eddy variables and in the limit $\epsilon_T\to0$.
This section considers the limit $\epsilon\to0$ while holding $\epsilon_T$ fixed.

The two approximate models are derived using standard approximations for slow-fast systems \citep{PK74,PS08,FW12}.
The presentation here follows the convenient review found in \cite{BGTV16}; the formulas are derived in a straightforward manner using formal asymptotic methods applied to the backwards Kolmogorov equation for the system (for details, see the appendices of \cite{BGTV16}).

The first approximation is derived via simple averaging.
In the limit $\epsilon\to0$ the eddy variables are well approximated as solutions to (\ref{eqn:v})--(\ref{eqn:S}) with $x$ and $y$ considered constant.
Curiously, although the full system (\ref{eqn:True}) has a smooth invariant measure the system (\ref{eqn:v})--(\ref{eqn:S}) does not: the long-time limiting distribution of $v$, $T$, and $S$ is jointly Gaussian with a singular covariance matrix.
In light of this, the following noise-augmented system is considered instead
\begin{subequations}\label{eqn:Fast}
\begin{align}
\text{d}v &=-\frac{v}{\epsilon}\text{d}v+\sqrt{\frac{2}{\epsilon}}\text{d}W_v\\
\text{d}T&=-\frac{1}{\epsilon}\left[T+2P^2vx\right]\text{d}t+\sqrt{\frac{2}{\epsilon}}\sigma_\epsilon\text{d}W_T\\
\text{d}S&=-\frac{1}{\epsilon}\left[S+2P^2vy\right]\text{d}t+\sqrt{\frac{2}{\epsilon}}\sigma_\epsilon\text{d}W_S
\end{align}
\end{subequations}
and the limit $\sigma_\epsilon\to0$ is taken after the fact.

The invariant measure of (\ref{eqn:Fast}) is Gaussian with zero mean and covariance
\begin{equation}\label{eqn:StatCov}
\left[\begin{array}{c c c}
1&-P^2x & -P^2y\\
-P^2x & 2P^4 x^2+\sigma_\epsilon^2 & 2P^4xy\\
-P^2y & 2P^4xy & 2P^4 y^2+\sigma_\epsilon^2
\end{array}\right].
\end{equation}
The averages of the terms $vT$ and $vS$ in the slow equations with respect to the invariant measure of the fast system are simply $-P^2x$ and $-P^2y$, respectively.
It is worth noting that these values are independent of the auxiliary noise amplitude $\sigma_\epsilon$.
Inserting these into the slow equations leads to the following approximate model\\
\begin{center}{\bf Deterministic Approximation}\end{center}
\begin{subequations}\label{eqn:Det}
\begin{align}
\text{d}x &= \left[-\frac{1}{\epsilon_T}(x-1)-[1+P_a(x-y)^2]x - 4P^2x\right]\text{d}t + \sqrt{\frac{1}{\epsilon_T}}\sigma_x\text{d}W_x\label{eqn:xDet}\\
\text{d}y &= \left[1-[1+P_a(x-y)^2]y - 4P^2y\right]\text{d}t +\sigma_y\text{d}W_y\label{eqn:yDet}.
\end{align}
\end{subequations}
The model (\ref{eqn:Det}) is referred to as the `deterministic' or `averaged' approximation since it models the eddy terms $vT$ and $vS$ as deterministic functions of $x$ and $y$.
It is straightforward to verify that this model is ergodic under the classical conditions of \cite{Khasminskii12}.\\

As described in \cite{BGTV16}, one can derive equations that approximate the variations of the true solution to (\ref{eqn:True}) around the solution of the approximate model (\ref{eqn:Det}).
Combining the equations for the variations with the deterministic approximation leads to further corrections in both the drift and diffusion, of order $\epsilon$ and $\sqrt{\epsilon}$, respectively.
The drift correction is significantly smaller than the leading-order drift.
But the leading-order diffusion terms in the $x$ and $y$ equations are of order $\approx0.1$, and corrections of order $\sqrt{\epsilon}$ may be of  comparable magnitude.

In order to compute the diffusion corrections, it is convenient to define some notation.
Let $\bm{Y}=(v,T,S)^T$ denote the solution to the noise-augmented system (\ref{eqn:Fast}).
Define constant matrices
\[\text{\bf M} = -\frac{1}{\epsilon}\left[\begin{array}{c c c}
1&0 &0\\
2P^2x&1&0\\
2P^2y&0&1
\end{array}\right],\;\;\text{\bf G} = \sqrt{\frac{2}{\epsilon}}\left[\begin{array}{c c c}
1&0 &0\\
0&\sigma_\epsilon&0\\
0&0&\sigma_\epsilon
\end{array}\right]\]
such that the fast system (\ref{eqn:Fast}) may be written d$\bm{Y} =$ {\bf M}$\bm{Y} + $ {\bf G}d$\bm{W}$.
The solution is thus
\begin{equation}\label{eqn:VOU}
\bm{Y}(\tau) = e^{\text{\bf M}\tau}\bm{Y}_0 + \int_0^\tau e^{\text{\bf M}(\tau-s)}\text{\bf G}\text{d}\bm{W}.
\end{equation}
The deviations of the eddy terms $vT$ and $vS$ from their conditional means are denoted
\[\bm{f}(x,y,\bm{Y}) = \left(\begin{array}{c}vT+4P^2x\\vS+4P^2y\end{array}\right).\]
According to \cite{BGTV16}, the diffusion-corrected model for the slow variables has the form
\begin{align*}
\text{d}x &= \left[-\frac{1}{\epsilon_T}(x-1)-[1+P_a(x-y)^2]x - 4P^2x\right]\text{d}t\\
&\qquad\qquad\qquad\qquad +\sqrt{\epsilon}a_{xx}(x,y)\text{d}\hat{W}_x+\sqrt{\epsilon}a_{xy}(x,y)\text{d}\hat{W}_y+ \sqrt{\frac{1}{\epsilon_T}}\sigma_x\text{d}W_x\\
\text{d}y &= \left[1-[1+P_a(x-y)^2]y - 4P^2y\right]\text{d}t \\
&\qquad\qquad\qquad\qquad +\sqrt{\epsilon}a_{yx}(x,y)\text{d}\hat{W}_x+\sqrt{\epsilon}a_{yy}(x,y)\text{d}\hat{W}_y+\sigma_y\text{d}W_y.
\end{align*}
where the matrix
\[\text{\bf A} = \left[\begin{array}{c c}a_{xx}&a_{xy}\\a_{yx}&a_{yy}\end{array}\right]\]
is any square root of the following symmetric positive definite matrix
\[\text{\bf C} = \int_0^\infty\mathbb{E}^{\bm{Y}_0}\left[\mathbb{E}^{\bm{Y}(\tau)}\left[\bm{f}(x,y,\bm{Y}(\tau))\bm{f}^T(x,y,\bm{Y}_0)+\bm{f}(x,y,\bm{Y}_0)\bm{f}^T(x,y,\bm{Y}(\tau))\right]\right]\text{d}\tau.\]
The matrix {\bf C} is the integral of the time-lagged auto-covariance of $\bm{f}$ with $x$ and $y$ considered constant.
In the above expression, $\mathbb{E}^{\bm{Y}(\tau)}$ denotes the expectation on $\bm{Y}(\tau)$ conditioned on the initial condition $Y_0$; the distribution is Gaussian with mean and covariance implied by (\ref{eqn:VOU}).
$\mathbb{E}^{\bm{Y}_0}$ denotes expectation on $\bm{Y}_0$ whose distribution is the stationary distribution of the fast process, in this case a zero-mean Gaussian with covariance (\ref{eqn:StatCov}).
The calculation for the system under consideration here is particularly straightforward since it requires only higher moments of jointly-Gaussian variables.
The matrix {\bf C} is found to have the form
\[\text{\bf C} = \left[\begin{array}{c c}
16(5P^4x^2+\sigma_\epsilon^2) & 80 P^4xy\\
80 P^4xy & 16(5P^4y^2+\sigma_\epsilon^2)
\end{array}\right].\]
In this case (unlike the leading-order drift term) the limit $\sigma_\epsilon\to0$ is singular in the sense that the matrix {\bf C} becomes positive semi-definite.
Nevertheless, a square root matrix {\bf A} exists; in the limit $\sigma_\epsilon\to0$ it has the form
\[\text{\bf A} = 4\sqrt{5}P^2\left[\begin{array}{c c}x&0\\y&0\end{array}\right].\]
The model for the slow variables with leading-order drift and diffusion corrections (but ignoring the order-$\epsilon$ drift correction) is thus\\
\begin{center}{\bf Gaussian Stochastic Approximation}\end{center}
\begin{subequations}\label{eqn:Stoch}
\begin{gather}
\text{d}x = -\left[\frac{1}{\epsilon_T}(x-1)+[1+P_a(x-y)^2]x + 4P^2x\right]\text{d}t+4\sqrt{5\epsilon}P^2x\text{d}\hat{W}+ \sqrt{\frac{1}{\epsilon_T}}\sigma_x\text{d}W_x\\
\text{d}y = \left[1-[1+P_a(x-y)^2]y - 4P^2y\right]\text{d}t +4\sqrt{5\epsilon}P^2y\text{d}\hat{W}+\sigma_y\text{d}W_y.
\end{gather}
\end{subequations}
For $x\approx1$ the noise amplitude associated with the eddies is $\approx0.16$, which is slightly larger than the `atmospheric' noise $\sigma_\epsilon/\sqrt{\epsilon_T}=0.1$.
The order-$\epsilon$ drift corrections have also been calculated, but they are small in comparison with the leading-order terms, and have been left out of the model for simplicity.
This system of SDEs is interpreted in the Ito sense; while the drift corrections in slow-fast systems with one slow degree of freedom can be interpreted as a correction from Stratonovich to Ito, this is no longer generally true in systems with multiple slow degrees of freedom \citep{PS08,FW12}.
It is straightforward to verify that this model is ergodic under the classical conditions of \cite{Khasminskii12}.

It is interesting to note that the Gaussian stochastic model replaces the eddy terms $vT$ and $vS$ by $-4P^2x($d$t + \sqrt{5\epsilon}$d$\hat{W})$ and $-4P^2y($d$t + \sqrt{5\epsilon}$d$\hat{W})$.
This form of subgrid-scale parameterization is qualitatively the same as that proposed in \cite{BMP99}, where it was proposed to multiply a deterministic parameterization (here $-4P^2x$) by a stochastic process (here $1+\sqrt{5\epsilon}\dot{W}$).
This style of stochastic parameterization has been widely used in atmospheric models (\cite{BernerEtAl16} provides a review), and much has been made of the role of multiplicative noise by, e.g.~\cite{SNPS05}.
The above derivation gives an example where this style of \emph{ad hoc} parameterization is rigorously justified, though multiplicative noise with a linear coefficient is certainly not the universal form of eddy-induced noise \citep[see e.g.][for a counterexample]{MC11}.

Recall that for the parameters (\ref{eqn:parameters}) the system (\ref{eqn:True}) has only three equilibria, two of which are stable.
The equilibria all have $v,T,S=0$, and the stable equilibria occur at $(x,y)\approx(0.989,0.22)$ and $(x,y)\approx(0.998,1.00)$.
The deterministic and Gaussian stochastic models have the same drift, which has only one equilibrium at $(x,y)\approx(0.974,0.093)$.
As will be verified by the results in \S\ref{sec:results}, the inclusion of nonlinear eddy effects completely changes the regime of the ocean model from a regime of multiple equilibria to a regime with a single stable equilibrium.

The averaged drift has a single stable equilibrium for all $P$ greater than approximately 0.117; below this value the drift undergoes a saddle-node bifurcation that creates a pair of equilibria near $x=1$ and $y=1$.
To achieve such small values of $P$ would require reducing the eddy velocity scale from 10 cm/s to 1 cm/s, which is unrealistically small.
The approximate models derived in this section show that the mean effect of eddies is linear and diffusive.
Since a linear diffusive effect is already present in the equations (the terms $-x$ and $-y$ in (\ref{eqn:x}) and (\ref{eqn:y})), the mean eddy effect could be viewed as a double-counting of eddy-induced diffusive exchange between the boxes.
This can be rectified by eliminating the mean diffusion terms, and such a model is formulated and studied in \S\ref{sec:NoMeanDiffusion}.
By avoiding a double-counting of diffusive exchange, the model in \S\ref{sec:NoMeanDiffusion} allows multiple equilibria with small, yet realistic eddy amplitudes.

\section{Numerical Methods}
\label{sec:methods}
Numerical methods are needed to compare the qualitative behavior of the three models (\ref{eqn:True}), (\ref{eqn:Det}), and (\ref{eqn:Stoch})	.
Many methods are derived based on the assumption that the drift is globally-Lipschitz \citep{KP92}, which is not the case here.
Several more recent investigations have analyzed numerical methods for SDEs whose drift satisfies a one-sided Lipschitz condition (e.g.~\citet{HMS02} and \citet{MS13}), but none of the models in consideration here satisfy such a condition.
A method appropriate to polynomial drifts is derived by \citet{LMS07}, but their analysis requires an invertible diffusion matrix, which the model (\ref{eqn:True}) does not have.
The Euler-Maruyama method may be appropriate, but is known to behave poorly in problems with polynomial drift \citep{MSH02,HJK11}.
In light of this, the `backward Euler' (BE) method is used here for all three models.
For a general system of SDEs of the form
\[d\bm{X} = \bm{b}(\bm{X})\text{d}t + \mathbf{\Sigma}(\bm{X})\text{d}\bm{W}\]
the BE method takes the following form
\begin{equation}\label{eqn:BE}
\bm{X}_{n+1} - \Delta t\bm{b}(\bm{X}_{n+1})= \bm{X}_n + \mathbf{\Sigma}(\bm{X}_n)\Delta\bm{W}_n
\end{equation}
where $\Delta t$ is the time step.
In every simulation presented here $\Delta t = 2\times10^{-6}$, which is significantly smaller than the smallest time scale of the system $\epsilon=2\times10^{-4}$.
\citet{MSH02} prove that the method is ergodic (for sufficiently small $\Delta t$) and that the invariant measure of the numerical method converges to that of the SDE as $\Delta t\to0$.
Though the analysis of \citet{MSH02} focuses on models with additive noise, the BE method is nevertheless applied here to the model (\ref{eqn:Stoch}) with multiplicative noise.

For the model (\ref{eqn:True}), a two-step process is used to generate solutions of the nonlinear system of equations (\ref{eqn:BE}).
First, an asymptotic approximation in the limit $\Delta t\to0$ is computed that has the form $\bm{X}_* = \bm{X}_n + \mathbf{\Sigma}(\bm{X}_n)\Delta\bm{W}_n + \mathcal{O}(\Delta t)$; this approximation is followed by a single Newton step.
For the systems (\ref{eqn:Det}) and (\ref{eqn:Stoch}), approximate solutions to the the nonlinear systems were generated using 10 fixed-point iterations started at $\bm{X}_n$.
Given the small step size, the resulting approximations solve their respective nonlinear systems with high accuracy; the residuals are typically on the order of $10^{-11}$.\\

\section{Results}
\label{sec:results}
\subsection{Climatology}
\begin{figure}[t]
  \centering
  \includegraphics[width=\textwidth]{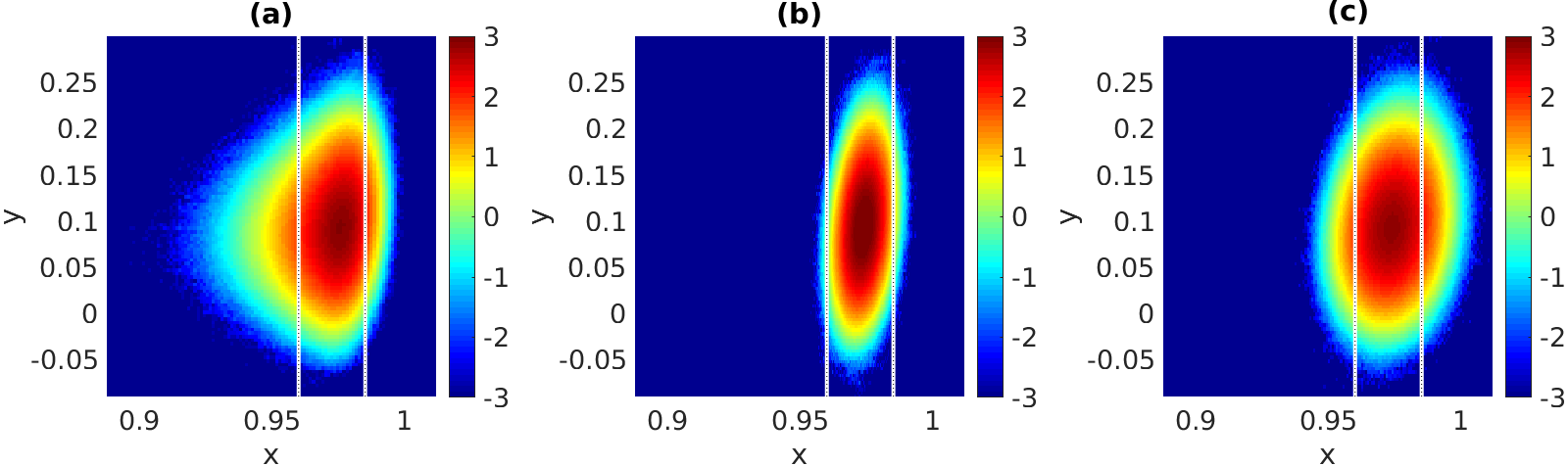}
  \caption{Base-10 logarithm of the climatological joint probability density functions of $x$ and $y$ for (a) Full model (\ref{eqn:True}), (b) Deterministic approximation (\ref{eqn:Det}), and (c) Gaussian-stochastic approximation (\ref{eqn:Stoch}). The vertical lines are placed at $x=0.96$ and $x=0.985$.}\label{fig:PDFs}
\end{figure}

A suite of 10,000 independent simulations was run starting from $x,\,y,\,v,\,T,\,S=0$ at $t=0$.
Data were saved for the time interval $t\in[4,10]$, saving every 100$^\text{th}$ time step for a spacing of $2\times10^{-4}$.
The mean and covariance appeared to have stabilized by $t=4$, suggesting that the data in $t\in[4,10]$ represents the stationary climatological distribution of the system.
Recalling that the dimensional time unit is 219 years, this amounts to 1,314 years of data saved approximately twice per month.
The three models all have the same mean of $(x,y)\approx(0.974,0.094)$, which is very close to the equilibrium of the deterministic and Gaussian stochastic models at $(0.974,0.093)$.
All three models have the same marginal standard deviation of $y$ approximately equal to $0.034$.
This can be explained by the fact that the amplitude of the eddy noise in the $y$ equation is estimated in the Gaussian stochastic model to be $4\sqrt{5\epsilon}P^2y\approx0.015$ for $y=0.093$, which is much less than the atmospheric noise with amplitude $\sigma_y=0.15$.

The climatological distributions of the models differ in other respects.
For example, the marginal standard deviation of $x$ is 0.0063, 0.0035, and 0.0065 in the full, deterministic, and Gaussian stochastic models, respectively.
The eddy noise in the $x$ equation is of comparable size to the atmospheric noise, and has a significant impact on the variability; the deterministic model lacks this eddy noise, and has too little variability.
The lack of eddy noise in the $x$ equation of the deterministic model also leads to an overestimate of the correlation between $x$ and $y$: the full and Gaussian stochastic models have correlations $0.15$ and $0.14$, respectively, while the deterministic model has correlation $0.23$.
The most-probable values of the distributions are $(x,y)\approx(0.976,0.092)$ for the full model, $(0.974,0.091)$ for the deterministic model, and $(0.973,0.093)$ for the Gaussian stochastic model; the differences in the $y$ value are negligible, but the differences in the $x$ value are up to half of a standard deviation.

Time-lagged correlation functions were computed, for example Corr$[x(t),x(t+\tau)] = C(\tau)$ (stationarity is assumed).
The correlation functions are all very similar across the models (not shown).
The correlation functions all decay monotonically to zero, so it is natural to define a decorrelation time by $\int_0^\infty C(\tau)\text{d}\tau$.
The correlation functions for $y$ in all three models are very similar, with decorrelation time approximately 22 years.
The correlation function for $x$ exhibits similar rapid initial decay in all three models.
The correlation function for $x$ in the deterministic model has a long tail, with larger long-lag correlations than the other two models, leading to a decorrelation time of 1.6 years, which is longer than the decorrelation times of the full model and Gaussian stochastic model, both of  which are approximately 1 year.

A simple binning procedure was used to generate approximations to the climatological probability density function (pdf) for each model; results are shown in Fig.~\ref{fig:PDFs}, with panels (a)--(c) presenting the full model, deterministic model, and Gaussian stochastic model, respectively.
It has already been noted that the three models have the same marginal variance for $y$, and indeed the range of $y$ in the three models is quite similar.
The deterministic model is clearly under-dispersed with respect to $x$. 
The climatological distribution of the Gaussian stochastic model has a more-accurate core, but is not skewed in the same way as the full model.

It is possible that minor deficiencies near the core of the distribution could be corrected by adding order-$\epsilon$ corrections to the drift of the Gaussian stochastic model, but the results of \cite{BGTV16} indicate that such corrections will not generate correct rare-event probabilities even in the limit $\epsilon\to0$.
To emphasize differences in the rare event probabilities, the probabilities of $x\le0.96$ and $x\ge0.985$ were calculated for the three models (these $x$ values are indicated by vertical lines in Fig.~\ref{fig:PDFs}).
The small-event probabilities are 0.039 for the full model, less than $10^{-4}$ for the deterministic model, and $0.022$ for the Gaussian stochastic model.
The large-event probabilities are 0.016 for the full model, less than $10^{-3}$ for the deterministic model, and $0.048$ for the Gaussian stochastic model.
Not surprisingly, the deterministic approximation has too-small rare event probabilities.
The Gaussian-stochastic model is more accurate, but is still incorrect by more than a factor of 2 in each case.

In order to verify that the system does not remain near the stable equilibrium of (\ref{eqn:True}) at $(x,y)\approx(1,1)$, a set of 1,000 simulations of (\ref{eqn:True}) was run with initial condition $(x,y,v,T,S)=(1,1,0,0,0)$.
These simulations were again run for the interval $t\in[0,10]$, saving the output from $t\in[4,10]$.
The stationary distribution did not display a secondary peak near $(1,1)$, indicating that the two stable equilibria of the full model are largely irrelevant to the dynamics of the system.\\

\subsection{Rare event forecasting}
\begin{figure}[t]
  \centering
  \includegraphics[width=\textwidth]{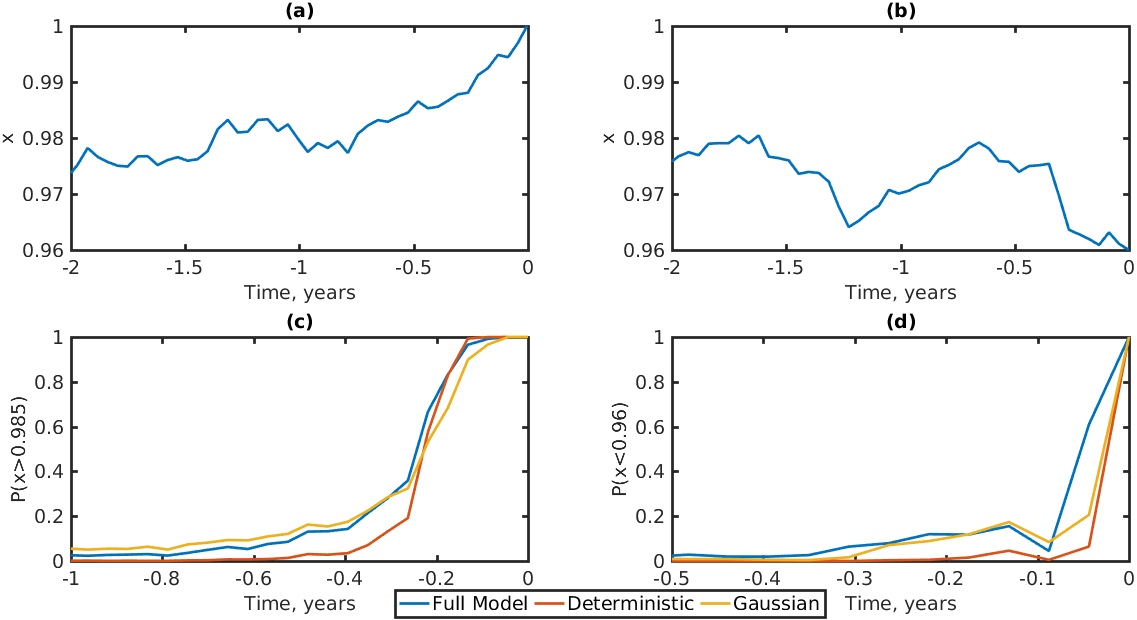}
  \caption{(a) and (b): $x$ trajectories of the full model (\ref{eqn:True}). (c) probability that $x>0.985$ at $t=0$, and (d) probability that $x<0.96$ at $t=0$ for forecasts initialized from the trajectories in (a) and (b), respectively. Note that the time axes in (a) and (b) are different from those in (c) and (d).}\label{fig:Forecasts}
\end{figure}

The previous section examined only the stationary climatological distributions of the three models.
Within a climate prediction scenario, short-term behavior is also important.
Given that the climatological distributions differ mainly in their rare event probabilities, a separate set of experiments was used to investigate the ability of the models to predict rare events over a shorter time interval.
The goal was to test how accurately the approximate models forecast the probability of the unusually large and small $x$ values over a range of forecast lead times.
Two trajectories of the system (\ref{eqn:True}) were selected out of the 10,000 discussed above: one with a value of $x\le0.96$ and one with a value of $x\ge0.985$.
These trajectories are shown in Fig.~\ref{fig:Forecasts} panels (a) and (b).
Note that the large-$x$ trajectory passes the threshold of $0.985$ approximately half a year before the final time, whereas the small-$x$ trajectory crosses the $0.96$ threshold only at the last time step.
For a single lead time, e.g.~2 years before the rare event, ensembles of 10,000 independent forecasts for all three models were initialized from the true trajectory.
At the time of the rare event these ensembles were used to estimate the probability of $x\le0.96$ (or $x\ge0.985$).
The probabilities shown in Fig.~\ref{fig:Forecasts}c correspond to the large-$x$ trajectory, and those in Fig.~\ref{fig:Forecasts}c correspond to the small-$x$ trajectory.
Since the large-$x$ trajectory crosses the threshold nearly half a year before the final time, all 10,000 of the forecasts initialized at any lead time less than half a year in advance are already above threshold; nevertheless, the probability at the final time is less than one because many of the trajectories cross the threshold back towards smaller values of $x$.

In both cases the forecast by the deterministic model is severely untrustworthy, assigning extremely low rare-event probability even when initialized close to or even beyond the threshold.
The rare-event probability forecast by the Gaussian stochastic model, in contrast, begins to increase from its climatological value at approximately the same time that the true forecast probability begins to increase, between 0.8 and 0.6 years in advance for the large-$x$ event and around 0.3 years in advance for the small-$x$ event.
Although the actual probability assigned by the Gaussian stochastic model at relatively long lead times is incorrect, the fact that it begins to increase at the right time could still be used qualitatively to predict whether the model is getting close to a rare event.
Once the probability of a rare event increases past about 20\%, the Gaussian stochastic model uniformly under-predicts the correct probability, despite having over-predicted the climatological probability for $x>0.985$.
For example, with a lead time of about 2.5 months the Gaussian stochastic model predicts the large-$x$ event with probability only 53\% while the true probability is in fact 67\%; with a lead time of half a month the Gaussian stochastic model predicts the small-$x$ event with probability only 21\% while the true probability is 61\%.

In summary, the deterministic model is essentially useless for rare-event forecasting, while the Gaussian stochastic model is only qualitatively useful, predicting whether a rare event is more likely but not with a robust uncertainty estimate.\\

\section{A model without mean diffusion}
\label{sec:NoMeanDiffusion}
As noted at the end of \S\ref{sec:Approx}, the averaged effect of the eddies is linear and diffusive.
Linear diffusive terms are already included in the budgets of heat and salt, with the result that the averaged models have only one stable equilibrium unless the eddies are assumed to be extremely weak, with velocities on the order of 1 cm/s.
If one assumes that linear diffusive exchange between the boxes is entirely eddy-driven then one can drop the mean diffusion terms from the governing equations of the full model, i.e.~equations (\ref{eqn:x}) and (\ref{eqn:y}) are changed to

\begin{equation}
\text{d}x = \left[-\frac{1}{\epsilon_T}(x-1)-P_a(x-y)^2x + 4vT\right]\text{d}t + \sqrt{\frac{1}{\epsilon_T}}\sigma_x\text{d}W_x\label{eqn:6x}
\end{equation}
and
\begin{equation}
\text{d}y = \left[1-P_a(x-y)^2y + 4vS\right]\text{d}t +\sigma_y\text{d}W_y\label{eqn:6y}
\end{equation}
respectively.
The eddy reductions proceed as before, so that the $-x$ and $-y$ terms are similarly dropped from the deterministic (\ref{eqn:Det}) and Gaussian (\ref{eqn:Stoch}) models.
The resulting model is much more amenable to multiple equilibria.
For $P$ greater than about 0.514 there is a single stable equilibrium with $x\approx1$ and $y\approx0.25$.
Below this value of $P$ the system undergoes a saddle-node bifurcation that creates a pair of equilibria near $(x,y)=(1,1)$; the saddle then moves down towards the original equilibrium, which it joins in a reverse saddle-node bifurcation at $P$ approximately 0.301, below which there remains only a single equilibrium.
We investigate the system at a value of $P_e=32$, i.e.~$P\approx0.45$, where there are three equilibria: a stable one at $(.99,.24)$, a saddle at $(1.00,.65)$, and another stable one at $(1.00,1.11)$.

\subsection{Ergodicity}
Recall that there are two conditions for ergodicity of hypoelliptic SDEs in \cite{MSH02}.
The first condition is that there is an inner-product norm $\|\bm{\cdot}\|$ such that $\langle\bm{u},\bm{F}(\bm{u})\rangle\le\alpha-\beta\|\bm{u}\|^2$ for some $\alpha,\beta>0$ where $\bm{u}$ is a vector containing the dependent variables and $\bm{F}(\bm{u})$ is the drift.
The second condition is that the vectors $\left\{\bm{\rho}_i,[[\bm{F},\bm{\rho}_j],\bm{\rho}_k]\right\}$ span $\mathbb{R}^5$ where $\bm{\rho}_i$, $i=1,\,2,\,3$ are the columns of the diffusion matrix, and $[\mathbf{\cdot},\mathbf{\cdot}]$ is a Lie bracket.
It is straightforward to verify that the second condition is met in this model in the same way that it is met in the original model (\ref{eqn:True}).

The first condition is more difficult.
We will use the inner product $\langle\bm{u},\bm{v}\rangle = u_1v_1+u_2v_2+\epsilon u_3v_3+(2\epsilon/P^2)(u_4v_4+u_5v_5)$, so we must show that there are $\alpha,\beta>0$ such that
\[\langle\bm{u},\bm{F}(\bm{u})\rangle -\alpha+\beta\|\bm{u}\|^2\le0\]
i.e.
\begin{multline*}-\alpha+ y-x(x-1)/\epsilon_T -P_a(x-y)^2(x^2+y^2) - v^2 -(2/P^2)(T^2+S^2)\\
+\beta(x^2+y^2+\epsilon v^2 + (2\epsilon/P^2)(T^2+S^2))\le0.
\end{multline*}
The terms involving the eddy variables ($v$, $T$, and $S$) will clearly pose no problem provided that $\beta < \epsilon^{-1}$.
It therefore remains to see whether one can choose $\alpha, \beta$ such that
\[-\alpha+ y-x(x-1)/\epsilon_T -P_a(x-y)^2(x^2+y^2)+\beta(x^2+y^2)\le0.\]
Consider the behavior along a line through the origin in the $(x,y)$ plane: along any line except $y=x$ the function is a quartic polynomial that can be made negative by choosing $\alpha$ sufficiently large.
Along the line $y=x$ the condition reduces to
\[-\alpha+ x -x(x-1)/\epsilon_T +2\beta x^2\le0.\]
As long as $\beta< \epsilon_T/2$ it will be possible to choose $\alpha$ sufficiently large that this condition is met.
The model without mean diffusion terms is therefore still ergodic.
Ergodicity is important because it implies that there is a single climatological distribution independent of the initial condition; the conditions of \cite{MSH02} further guarantee that the distribution collapses exponentially quickly towards the climatological distribution.

\subsection{Numerical experiments}
\begin{figure}[t]
  \centering
  \includegraphics[width=\textwidth]{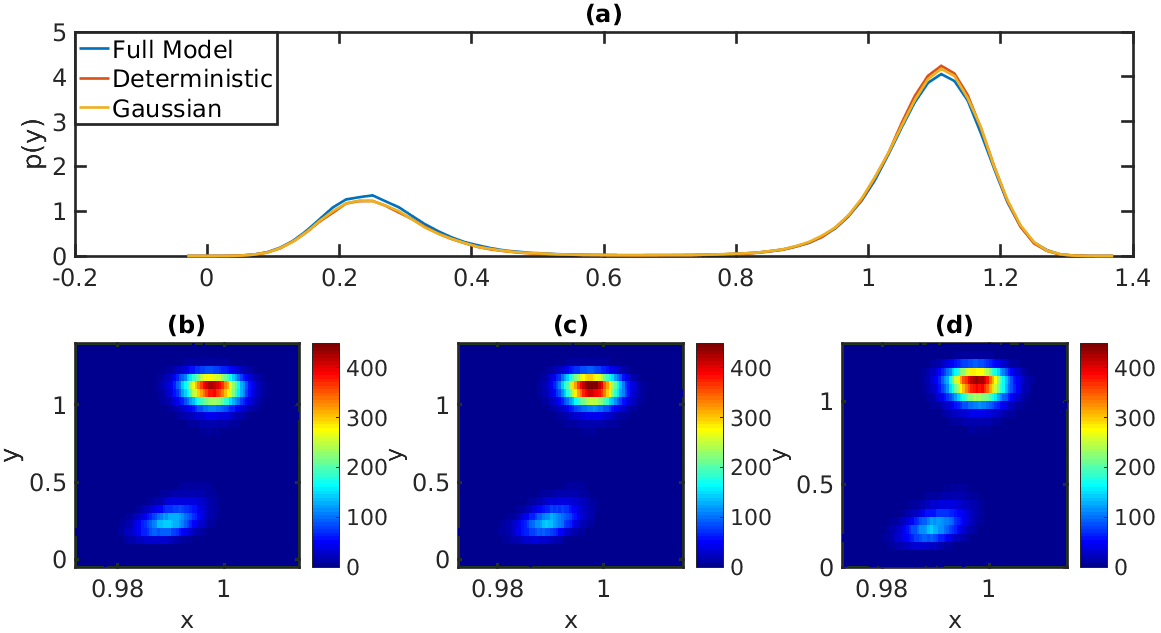}
  \caption{(a) Climatological marginal probability density functions $p(y)$ for the three models without mean diffusion. Climatological joint probability density functions $p(x,y)$ for (b) Full model, (c) Deterministic approximation, and (d) Gaussian-stochastic approximation.}\label{fig:NoDiff}
\end{figure}

Ensemble simulations for the three models without mean diffusion were run with 1000 ensemble members each.
The deterministic and Gaussian approximate models were initialized with $x=1$, $y=0.6$, while the full model was initialized with $x=1$, $y=0.65$, and $v,T,S=0$.
After a burn-in of 4 nondimensional time units, the simulations were run for 500 more time units, i.e.~about 110,000 years.
Although the models are geometrically ergodic, with distributions collapsing exponentially quickly towards the invariant distribution, this was not enough time for the approximate models to reach the invariant distribution.
These models were then extended for a further 500 time units, during which time their distributions converged.
The full model was initialized closer to the saddle point, so its distribution converged within the first 504 time units.

The climatological distributions of the slow variables are shown in Fig.~\ref{fig:NoDiff}.
Panel (a) shows the marginal $y$ distributions of the three models, while panels (b)--(d) show the joint $(x,y)$ distributions.
The three models are remarkably similar.
Though $\epsilon$ is the same as in the previous case, $P$ is smaller.
The diffusion correction in the $x$ equation of the Gaussian-stochastic approximation has amplitude $4\sqrt{5\epsilon}P^2x\approx0.026$ which is smaller than the atmospheric noise amplitude $\sigma_x/\sqrt{\epsilon_T} = 0.1$; the atmospheric noise similarly dominates the $y$ equation.
As a result, the effects of eddy noise are not seen in the equilibrium distributions of the three models.\\

\begin{figure}[t]
  \centering
  \includegraphics[width=\textwidth]{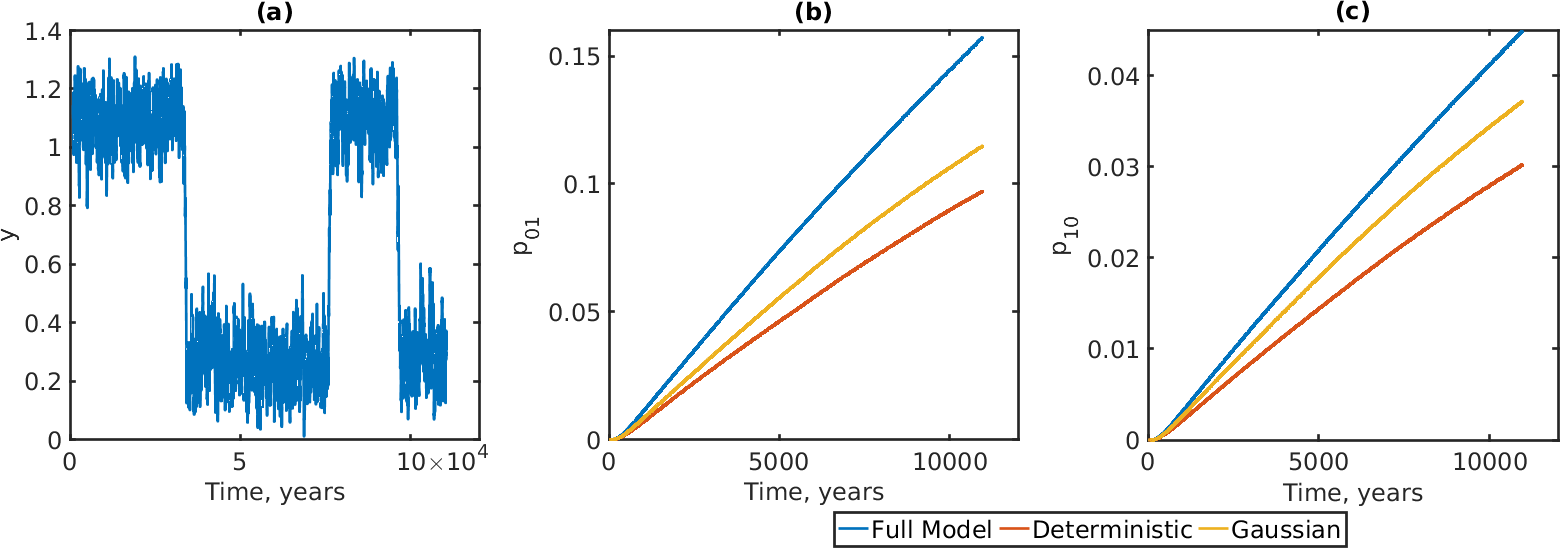}
  \caption{Regime transitions for the three models without mean diffusion. (a) A single $y(t)$ trajectory from the full system showing jumps between regimes. (b) The probability $p_{01}(\tau)$ of a transition from $y(t)<0.5$ to $y(t+\tau)>0.8$. (c) The probability $p_{10}(\tau)$ of a transition from $y(t)>0.8$ to $y(t+\tau)<0.5$.}\label{fig:TransProbs}
\end{figure}

The noise levels are low enough that the system trajectories make rare transitions between the neighborhoods of the two stable equilibria; Fig.~\ref{fig:TransProbs} panel (a) shows a system trajectory $y$ from the full model that jumps between regimes.
The rates and paths of these transitions are the subject of large deviation theory \citep{FW12}.
The methods of \citet{BGTV16} to analyze the transitions do not seem to apply directly here because of the inclusion of noise forcing in the slow dynamics.
In any case, it is not difficult to estimate the transition probabilities from simulations.
For practical purposes it was convenient to estimate the following probabilities $p_{01}(\tau) = P(y(t+\tau)>0.8\,| \,y(t) < 0.5)$ and $p_{10}(\tau) = P(y(t+\tau)<0.5\,|\,y(t) > 0.8)$.
These transition probabilities are plotted for the three models in panels (b) and (c), respectively.
The effects of differences in the eddy noise are clear: the deterministic model has the lowest transition probabilities; the Gaussian stochastic model has higher transition probabilities; the full model has the highest transition probabilities.

\section{Conclusions}
\label{sec:conclusions}
This paper formulates a stochastic two-box ocean model modeled after Stommel's (\citeyear{Stommel61}); the model consists of a system of 5 SDEs (\ref{eqn:True}). 
Previous stochastic Stommel models \citep[e.g.][]{Cessi94,VACTH01,MTL02,Monahan02,MC11}, modeled the atmospheric heat and freshwater forcing as Gaussian stochastic processes, and the exchange of heat and salt between the boxes as a nonlinear drift term corresponding to the large-scale overturning thermohaline circulation.
The novelty of the formulation here is that a fast, eddy-driven component is added to the the exchange between the boxes.
The terms modeling the eddy-driven exchange are quadratic products of approximately Gaussian random variables; products of jointly-Gaussian random fields were recently found to be an accurate model of eddy-driven exchanges in \cite{Grooms16}.

In more complete and complex ocean models, fast eddy effects are frequently modeled deterministically.
Stochastic parameterizations have recently been developed that multiply these deterministic eddy parameterizations by Gaussian random fields \citep{ACJPWZ16,JPZ17}, which is a popular approach for atmospheric models based on the work of \citet{BMP99} and \citet{SNPS05}.
Using methods of averaging and homogenization for slow-fast systems \citep{PS08,FW12,BGTV16}, two models were derived approximating the evolution of the slow components (the difference in heat and salt content of the two boxes).
The first model (\ref{eqn:Det}) replaces the fast eddy-driven exchange terms by a fixed `deterministic' drift term, analogous to the standard approach of deterministic parameterization in more complex ocean models.
The second model (\ref{eqn:Stoch}) adds an additional multiplicative noise term accounting for fast variations in the eddy-driven flux.
A suite of simulations of each of the three models was used to compare their qualitative behavior.
Two regimes were considered: one with a single stable equilibrium, and one with two stable equilibria.

The main results are as follows.
There is little qualitative difference in the core of the stationary distributions of the full, non-Gaussian model and the Gaussian multiplicative approximation.
In the regime with a single equilibrium the deterministic model has too little variability, but the Gaussian model gives an accurate climatological mean and covariance.
In the regime with two stable equilibria the climatological distribution of the three models is nearly the same.
In the regime with two stable equilibria the amplitude of the eddies is smaller than in the regime with a single equilibrium, which could perhaps account for the fact that the deterministic model is more accurate in the former regime.
Observational estimates suggest that up to 30\% of the variability of the Atlantic Meridional Overturning Circulation (AMOC) is driven by ocean eddies, with the rest driven by atmospheric noise \cite{HirschiEtAl13,SonnewaldEtAl13}.

Though the Gaussian stochastic model gives a good approximation of the core of the climatological distribution, the rare event probabilities are inaccurate.
In the single-equilibrium regime there is no clear trend in the behavior: the Gaussian model overestimates rare event probabilities on one side of the mean, and underestimates on the other side.
This inaccuracy manifests for short time, transient behavior too: even with a short lead time, the Gaussian model gives inaccurate predictions of the probability of a rare event.

In the regime with two stable equilibria the rare events of interest are the transitions between the two.
Despite the reduced amplitude of the eddies in this regime, clear differences were observed in the rates of transition from the neighborhood of one equilibrium to another: the deterministic model had the rarest transitions, and the Gaussian model still made transitions less frequently than the full model.

The goal of the investigation was to investigate the qualitative impacts of non-Gaussian eddy noise of the type observed by \citet{Grooms16} in a simple model, and to compare to models with Gaussian noise and without eddy noise.
The results suggest that Gaussian stochastic parameterizations may be able to successfully produce the day-to-day variability associated with the core of the climatological distribution, but that more accurate non-Gaussian models may be needed to correctly model rare events.
Such rare events include extreme behavior like droughts and heat waves, as well as abrupt transitions between climate regimes.
The impact of stochastic parameterizations on rare event distributions in climate models has only recently begun to be investigated \cite{TagleEtAl16}.

\section*{Acknowledgments}
IG is grateful to F.~Bouchet for a discussion of the methods from \cite{BGTV16} that are applied here. WB is supported by NSF EXTREEMS grant 1407340.

\end{document}